\documentclass[a4paper]{scrartcl}
       
\usepackage[utf8]{inputenc}
\usepackage[T1]{fontenc}
\usepackage[ngerman]{babel} 
\usepackage{lmodern}
\usepackage[onehalfspacing]{setspace}
\usepackage{natbib}
\usepackage{amsmath,amssymb,amsthm,amsfonts,amsbsy,latexsym}
\usepackage{graphicx}
\usepackage{geometry}
\usepackage[pdfpagelabels=true]{hyperref}

\newcommand{\E}{\mbox{I\negthinspace E}}
\newcommand{\B}{\mbox{I\negthinspace B}}

\newcommand{\R}{\mathbb{R}}

\DeclareMathOperator*{\argmin}{argmin}
\DeclareMathOperator*{\sign}{sign}

\providecommand{\keywords}[1]
{
  \small	
  \textbf{\textit{Keywords---}} #1
}

\numberwithin{equation}{section}

\begin{document}

\newtheorem{thm}{Theorem}[section]
\newtheorem{Def}{Definition}[section]
\newtheorem{lem}{Lemma}[section]
\newtheorem{rem}{Remark}[section]
\newtheorem{cor}{Corollary}[section]
\newtheorem{ex}{Example}[section]
\newtheorem*{bew}{Proof}

\title{Asymptotic linear expansion of regularized M-estimators}
\author{Tino Werner\footnote{Institute for Mathematics, Carl von Ossietzky University Oldenburg, P/O Box 2503, 26111 Oldenburg (Oldb), Germany}}
\maketitle

\begin{footnotesize} 

\begin{abstract}

Parametric high-dimensional regression analysis requires the usage of regularization terms to get interpretable models. The respective estimators can be regarded as regularized M-functionals which are naturally highly nonlinear. We study under which conditions these M-functionals are compactly differentiable, so that the corresponding estimators admit an asymptotically linear expansion. In a one-step construction, for a suitably consistent starting estimator, this linearization replaces solving optimization problems by evaluating the corresponding influence curves at the given data points. We show under which conditions the asymptotic linear expansion is valid and provide concrete examples of machine learning algorithms that fit into this framework. \end{abstract} 

\keywords{Asymptotic linear expansion, Regularized M-estimators, Influence curves} \ \\

\end{footnotesize}

\begin{small}

\section{Introduction} \

In the mid nineties, Robert Tibshirani succeeded in combining two important paradigms of fitting linear regression models, namely variable selection and shrinkage of the coefficients, in one single optimization problem, calling it the Lasso (Least Absolute Shrinkage and Selection Operator, \cite{tibsh96}). While already being superior to the former state-of-the-art procedures of Ridge regression and subset selection in terms of interpretability of the model and prediction accuracy (\cite{tibsh96}), its popularity grew when Efron et. al. embedded the Lasso into the framework of forward stagewise regression and provided the LARS algorithm (\cite{efron04}) which turned out to be more efficient than the implementations of Tibshirani and \cite{osborne00}. A very competetive algorithm has been developed by Friedman et al. in 2007, relying on the fact that despite the Lasso for multiple regression does not have a closed form solution, a simple Lasso just concerning one single predictor has. Therefore, they apply the so-called ''shooting'' algorithm to the Lasso and other suitable problems, which means that one repeatedly cycles through the variables, letting all others fixed to their values of the previous iteration and fits the partial residual, i.e., a coordinate-wise optimization is done (\cite{friedman07}). \ \\

Also going back to the nineties, the idea of combining several simple or ''weak'' models to get a ''strong'' model has been proposed. The first popular so-called ''Boosting'' algorithm was \texttt{AdaBoost} for solving classification problems (\cite{freund97}). In the path-breaking work of Friedman in 2001 (\cite{friedman01}), a gradient descent Boosting algorithm has been presented, where the negative gradient vector of the loss function, evaluated in the current model, is fitted in each iteration by a weak model. At the end of each iteration, the model is additively updated. A special attention was later given to the $L_2-$Boosting, that is the gradient Boosting with the squared loss, in seminal works of Bühlmann et al. (\cite{bu03}, \cite{bu06}, \cite{bu07}) where it is shown that it can be achieved by iteratively computing simple least squares fits. It turns out that $L_2-$Boosting (and other types of Boosting as well) performs variable selection and shrinks the coefficients, therefore it is strongly related to the Lasso, although the striking difference that the Lasso can at most select $n$ predictors if $p>n$ makes clear that they are not equal, as convincingly pointed out in \cite[Sec. 4.3]{bu06}.  \ \\

The advantages of both methods are clear: Instead of directly solving an optimization problem w.r.t. a $p-$dimensional parameter, one splits this task up into several simple optimization problems with a known closed-form solution, leading to a very high computational efficiency. More precisely, from the fact that $L_2-$Boosting performs variable selection, it can be implicitely seen as minimizing a regularized $L_2-$criterion. The variable selection property of Boosting is inherited from the evaluation of different simple models in the given loss function with preference for the best performing one. In fact, an aggregation over the columns is made, both in the fitting step as well as in the updating step. This may be regarded as a kind of ''measure'' on the columns of the data set which has been introduced and studied in the Ph.D. thesis of the author (\cite{TWphd}). \ \\ 

This contribution should go in the direction of a linearization of certain regularized M-estimators which corresponds to an aggregation over the rows of the data matrix. Such a linearization can be identified with a first-order expansion, invoking functional gradients. In fact, such functional gradients have already been mentioned in \cite[Sec. 2.1.1]{bu07} where the steepest gradient descent algorithm is formulated from a functional point of view. The gradient vector equals the G\^{a}teaux derivative of the empirical risk functional. Without spelling it out, Bühlmann exactly used the definition of the influence curve of a corresponding M-functional (\cite{hampel}), so the gradient vectors $U_i$ correspond to the influence curve evaluated at the data points $(X_i,Y_i)$. Those influence curves play a major role in robust statistics since they quantify the infinitesimal influence of a single observation to the estimator. Moreover, since one of the most important contributions of robust statistics is the identification of estimators with functionals  (\cite{huber}, \cite{hampel}, \cite{rieder}, \cite{maronna}), the influence curves which are nothing but functional first gradients linearize an estimator under certain conditions (\cite{mises}, \cite{reeds}, \cite{rieder}), up to a remainder term of order $o(n^{-1/2})$. Such an estimator is referred to as asymptotically linear estimator (ALE). \ \\

Therefore, the main question for this work is the following one: Can we make assumptions under which regularized M-estimators are asymptotically linear? And if we can, is it possible to embed the Lasso into this framework? \ \\

The G\^{a}teaux derivatives of certain penalized M-functionals have already been computed (\cite{alfons13}, \cite{alfons14}, \cite{medina}), but working with them requires to prove that the estimator can be linearized, i.e., that the remainder can be controlled uniformly under suitable distributional neighborhoods, which was only mentioned as an idea in terms of Fréchet differentiability (\cite{averbukh}) in \cite{medina} and which has not been done in the other references. There already exist results where the asymptotic linearity of penalized M-estimators has been shown (\cite{geer16}, \cite{laan15}) and even remarkable results concerning penalized M-estimators in the nonconvex setting (\cite{loh15}), but to the best of our knowledge, in particular, the case of a general theory for the asymptotic linearity of penalized M-estimators has not been covered by the aforementioned results so far.\ \\ 

In literature, robustness results by proving the boundedness of influence curves have been established in the special case of regularized kernel-based regression problems (\cite{christ04}). Since they require Fréchet-differentiability of the loss function which is not true for example for the $\epsilon-$insensitive loss, \cite{christ08} introduced Bouligand influences curves to show the robustness of support vector regression estimators. The asymptotic normality of kernel-based regression methods has been shown in \cite{hable12} who also used the framework of compact differentiability of the corresponding functional. \ \\

The rest of this article is organized as follows. In section \ref{prelim}, the general definition of $\mathcal{R}-$differentiability and relevant tools from robust statistics are introduced. In section \ref{compdiffsec}, results from \cite{rieder} are recapitulated which will be essential for the rest of this main section. Then we transfer the results for M-estimators to the case of regularized M-estimators. Subsection \ref{alesubs} is the main theoretical part where we show under which conditions a compactness assumption of the parameter space is reasonable and when the asymptotic linearity is valid. For non-differentiable penalty terms, we will heavily rely on an approximation lemma of \cite{medina}. Section \ref{exsec} is devoted to concrete examples. \ \\

\section{Preliminaries}\label{prelim} \

This section compiles the concepts needed for the main section. We recur to the abstract definition of $\mathcal{R}-$differentiability of maps between normed vector spaces. The second part contains the most important definitions of quantitative robust statistics like the influence curve, the $L_2-$differentiability of a parametric statistical model and asymptotically linear estimators. \ \\

\subsection{Functional derivatives} \

We begin with the definition of $\mathcal{R}-$differentiability of maps between (possibly infinite-dimensional) normed real vector spaces (cf. \cite{rieder}, \cite{averbukh}). 

\begin{Def}\label{rderiv} Let $X$, $Y$ be normed real vector spaces. A map $T: X \rightarrow Y$ is \textbf{$\mathcal{R}-$differenti-able}\index{Functional differentiability!$\mathcal{R}-$} in $x \in X$ if there exists $d_{\mathcal{R}}T(x) \in L(X,Y)$, $s_0>0$, such that \begin{center} $ \displaystyle T(x+sh)=T(x)+d_{\mathcal{R}}T(x)sh+\rho(sh) \ \forall |s| \le s_0 $, \end{center} where the remainder term $\rho$ satisfies the following conditions: 

\textbf{i)} $\rho(0)=0$, 

\textbf{ii)} $\rho \in \mathcal{R}(X,Y)$ where $\mathcal{R}(X,Y)$ is a real vector subspace of $Y^X$ with $\mathcal{R}(X,Y) \cap \mathcal{L}(X,Y)=\{0\}$. \ \\

Then the continuous, linear map (w.r.t. $h$ per definition) $d_{\mathcal{R}}T(x)$ is referred to as the \textbf{$\mathcal{R}-$derivative \index{Functional differentiability!$\mathcal{R}-$derivative} of $T$ at $x$}. \end{Def} \ 

The following definition taken from \cite[Sec. 1]{rieder} is helpful to distinguish between different types of $\mathcal{R}-$differentiation.  

\begin{Def}\label{rderivunif} Let $X$, $Y$, $T$ be as in definition \ref{rderiv}. Let $\mathcal{S}$ be a covering of $X$. Define the \textbf{remainder class} \begin{center} $ \displaystyle \mathcal{R}_{\mathcal{S}}(X,Y):=\left\{ \rho: X \rightarrow Y \ \bigg| \ \lim_{t \rightarrow 0}\left(\sup_{h \in S}\left(\frac{||\rho(th)||}{t} \right) \right)=\rho(0)=0 \ \forall S \in \mathcal{S} \right\}. $ \end{center} \end{Def} 

By definiton \ref{rderivunif}, it is clear that $\mathcal{R}-$differentiability can be seen as a linear approximation of some functional $T$ such that the remainder term converges uniformly on all sets $S \in \mathcal{S}$. Following, we now define three special concepts of $\mathcal{R}-$differentiability (cf. \cite{rieder}). We say that the functional $T$ is G\^{a}teaux or weakly differentiable resp. Hadamard or compactly differentiable resp. Fréchet or boundedly differentiable if the cover of $X$ consists of finite resp. compact resp. bounded sets. Trivially, bounded differentiability implies compact differentiability which implies G\^{a}teaux differentiability. The derivatives coincide in this case. Moreover, continuous G\^{a}teaux differentiability implies bounded differentiability. \ \\

For the application of an infinite-dimensional delta method (\cite[Thm. 1.3.3]{rieder}), it is necessary to investigate whether the \textbf{chain rule} holds for functional derivatives. The following theorem ( \cite[Prop. 1.2.6+Thm. 1.2.9]{rieder}, \cite[Thm. 1.6]{averbukh}) shows that for linear maps as approximations, the chain rule holds if and only if at least Hadamard differentiability holds.

\begin{thm}[\textbf{Chain rule}]\label{compchain} Let $X$, $Y$, $Z$ be normed real vector spaces and let $T: X \rightarrow Y$, $U: Y \rightarrow Z$. If $T$ and $U$ are compactly differentiable, then the chain rule holds, i.e., \begin{center} $ \displaystyle  d_H(U \circ T)(x)=d_HU(T(x)) \circ d_HT(x). $ \end{center} Conversely, if the chain rule holds, the maps are already compactly differentiable. \end{thm} \

The chain rule does not hold for G\^{a}teaux differentiable maps in general. Counterexamples can be found in \cite{averbukh} or \cite{frechet}. Thus, from the above mentioned concepts of functional derivatives one can state that compact differentiability is the weakest form of $\mathcal{R}-$differentiability such that the chain rule holds. \ \\

Of course, there exist examples where Hadamard-differentiability fails. One typical example are L-statistics where the underlying distribution has an unbounded support as pointed out in \cite{vaart98}. Then, compact differentiability is impossible w.r.t. $||\cdot ||_{\infty}$. Such functionals are written in the form (cf. \cite{beutner10}) \begin{center} $ \displaystyle T_g(F):=-\int xdg(F(x)) $, \end{center} so it is required that $g$ has a compact support in $]0,1[$ to ensure compact differentiability. It is shown in \cite{beutner16} that if the support of $g$ contains at least one of the boundary points of $[0,1]$, even the negative expectation value (by setting $g:=id$) is not compactly differentiable. The functional $T_g$ covers relevant statistical functionals like the Value at Risk or the Average Value at Risk as pointed out in \cite{beutner10}. In fact, \cite{kraetsch} stated that tail-dependent functionals are in general not compactly differentiable w.r.t. uniform norms. \ \\

\subsection{Basic concepts of quantitative robustness} \

Every real data analysis requires model assumptions. However, these assumptions are in general not fulfilled, hence the real data differ from data that would have been generated by the ideal model. Therefore, fitting models by using the real data can be seen as if one analyzes a contaminated data set which affects the quality of the fitted model. It is not desirable to exclude potential ''outliers'' from the data set (cf. \cite{hampel}) but to find strategies that downweight them, like iteratively reweighted least squares (IRWLS) (\cite{huber}). We now define the influence function (cf. \cite{hampel74}). \ 

\begin{Def}\label{IC} Let $X$ be a normed function space and let $\Theta$ be a normed real vector space. Let $T: X \rightarrow \Theta$ be a statistical functional. The \textbf{influence function or influence curve} of $T$ at $x$ for a probability measure $P$ is defined as the derivative \begin{center} $ \displaystyle IC(x,T,P):=\lim_{t \rightarrow 0}\left( \frac{T((1-t)P+t \delta_x)-T(P)}{t} \right)=\partial_t \left[T((1-t)P+t \delta_x) \right] \bigg|_{t=0} $ \end{center} where $\delta_x$ denotes the Dirac measure at $x$. \end{Def}  \

So, the influence curve is just a special G\^{a}teaux derivative with $h:=\delta_x-P$. The influence curve can be seen as an estimate for the infinitesimal influence of a single observation on the estimator. If the IC is unbounded, then a single observation can have an infinite impact on the resulting estimator which is of course not desirable. For robustness properties, it is necessary that the influence curve is at least bounded. In that case, the estimator is sometimes called B-robust (\cite{vaart98}). \ \\ 

The robustification of an estimator can be done by robustifying its influence function. Minimax results for optimal-robust influence curves have been established in several works (\cite{rieder}, \cite{rieder08}, \cite{hampel}, \cite{fraiman}). However, for guaranteeing optimality of these approaches, it is crucial that the estimator is asymptotically linear (\cite[Def. 4.2.16]{rieder}) and that the model be smooth enough, i.e., that it be $L_2-$differentiable (\cite{lecam1970}). \

\begin{Def}\label{l2diff} Let $\mathcal{P}:=\{P_{\theta} \ | \ \theta \in \Theta\}$ be a family of probability measures on some measurable space $(\Omega, \mathcal{A})$ and let $\Theta$ be a subset of $\R^p$. Then $\mathcal{P}$ is \textbf{$L_2-$differentiable at $\theta_0$}\index{$L_2-$differentiability} if there exists $\Lambda_{\theta_0} \in L_2^p(P_{\theta_0})$ such that \begin{center} $ \displaystyle  \left| \left| \sqrt{dP_{\theta_0+h}}-\sqrt{dP_{\theta_0}}\left(1+ \frac{1}{2} h^T \Lambda_{\theta_0} \right) \right| \right|_{L_2}=o(||h||) $ \end{center} for $h \rightarrow 0$. In this case, the function $\Lambda_{\theta_0}$ is the $L_2-$derivative and \begin{center} $ \displaystyle I_{\theta_0}:=\E_{\theta_0}[\Lambda_{\theta_0}\Lambda_{\theta_0}^T] $ \end{center} is the Fisher information of $\mathcal{P}$ at $\theta_0$. \end{Def} \

Note that the $L_2-$differentiability is a special case of the wider concept of $L_r-$differentiability (\cite{rieder01}). The $L_2-$differentiability holds for many distribution families, including normal location and scale families, Poisson families, Gamma families, and even for ARMA, ARCH and GPD families (\cite{rieder08}, \cite{pupa15}). A standard example of a distribution family that is not $L_2-$differentiable is the model $\mathcal{P}:=\{U([0,\theta]) \ | \ \theta \in \Theta\}$.  \ 

\begin{Def}\label{ale} Let $(\Omega^n, \mathcal{A}^n)$ be a measurable space and let $S_n: (\Omega^n, \mathcal{A}^n) \rightarrow (\R^p, \B^p)$ be an estimator. Then the sequence $(S_n)_n$ is \textbf{asymptotically linear at $P_{\theta_0}$}\index{Asymptotic linearity|(} if there exists an influence curve\index{Influence curve/function} $\psi_{\theta_0} \in \Psi_2(\theta_0)$ such that the expansion \begin{center} $ \displaystyle  S_n=\theta_0+\frac{1}{n} \sum_{i=1}^n \psi_{\theta_0}(x_i)+o_{P_{\theta_0}^n}(n^{-1/2}) $ \end{center} holds. The family $\Psi_2(\theta_0)$ of influence curves is defined by the set of all maps $\eta_{\theta_0}$ that satisfy the conditions \ \\

 \textbf{i)} $\eta_{\theta_0} \in L_2(P_{\theta_0})$, \ \ \ \textbf{ii)} $\E_{\theta_0}[\eta_{\theta_0}]=0$, \ \ \ \textbf{iii)} $\E_{\theta_0}[\eta_{\theta_0} \Lambda_{\theta_0}^T]=I_p $ \ \\  

where $I_p$ denotes the identity matrix of dimension $p \times p$. \end{Def} \ 

In this definition, condition i) is vital for integrability and for the application of a central limit theorem to conclude that $S_n$ is asymptotically normal, i.e., \vspace{0.1cm} \begin{center} $ \displaystyle \sqrt{n}(S_n-\theta_0) \circ P_{\theta_0}^n=\left(\frac{1}{\sqrt{n}} \sum_{i=1}^n \psi_{\theta_0}(x_i)+o_{P_{\theta_0}^n}(n^0) \right) \circ P_{\theta_0}^n \overset{w}{\longrightarrow} \mathcal{N}_p(0,\E_{\theta_0}[\psi_{\theta_0} \psi_{\theta_0}^T]) .$ \end{center}  \vspace{0.1cm} Condition ii) ensures unbiasedness of the asymptotically linear estimator.  The third condition leads to uniform unbiasedness (w.r.t. $\theta_0$), more precisely, if $\psi_{\theta_0}$ satisfies i) and ii), \cite[Lemma 4.2.18]{rieder} shows that the condition iii) is equivalent to \begin{center} $ \displaystyle  \sqrt{n}(S_n-\theta_0)(P_{\theta_0+t_n/\sqrt{n}}^n) \overset{w}{\longrightarrow} \mathcal{N}_p(t,\E_{\theta_0}[\psi_{\theta_0}\psi_{\theta_0}^T]) $ \end{center} for all $t_n \rightarrow t$ where $t_n$, $t \in \R^p$, so the asymptotic normality granted by a central limit theorem will hold locally uniformly over compacts (\cite{rieder}). \ \\

An extension of this concept arises if one wants to estimate the transformed parameter $\tau(\theta)$ leading to so-called ''partial'' influence curves in the terminology of \cite[Def. 4.2.10]{rieder}, \cite{rieder08}. \

\begin{Def} \label{alepic} Let $(\Omega^n, \mathcal{A}^n)$ be a measurable space and let $S_n: (\Omega^n, \mathcal{A}^n) \rightarrow (\R^q, \B^q)$ be an estimator for the transformed quantity of interest $\tau(\theta)$. Assume that $\tau: \Theta \rightarrow \R^q$ is differentiable at $\theta_0 \in \Theta$ where $\Theta \subset \R^p$ and $q \le p$. Denote the Jacobian by $\partial_{\theta_0} \tau=:D_{\theta_0} \in \R^{q \times p}$. Then the set of \textbf{partial influence curves (pIC)}\index{Influence curve/function!Partial} is defined by \begin{center} $ \displaystyle \Psi_2^D(\theta_0):=\{\eta_{\theta_0} \in L_2^q(P_{\theta_0}) \ | \ \E_{\theta_0}[\eta_{\theta_0}]=0, \  \E_{\theta_0}[\eta_{\theta_0} \Lambda_{\theta_0}^T]=D_{\theta_0}\} $. \end{center} Then the sequence $(S_n)_n$ is asymptotically linear\index{Asymptotic linearity|)} at $P_{\theta_0}$ if there exists a partial influence curve $\eta_{\theta_0} \in \Psi_2^D(\theta_0)$ such that the expansion \begin{center} $ \displaystyle  S_n=\tau(\theta_0)+\frac{1}{n} \sum_{i=1}^n \eta_{\theta_0}(x_i)+o_{P_{\theta_0}^n}(n^{-1/2}) $ \end{center} is valid. \end{Def} \ 

Since it holds that  \begin{center} $ \displaystyle \Psi_2^D(\theta_0)=\{D_{\theta_0} \psi_{\theta_0} \ | \ \psi_{\theta_0} \in \Psi_2(P_{\theta_0})\} $ \end{center} (see \cite[Rem. 4.2.11 e)]{rieder}), the asymptotically linear expansion of transformed estimators in terms of partial influence curves clearly mimicks the traditional delta-method. \ \\ 

Asymptotic linearity has been proven for example for asymptotically normal $M, R$ and $MD$ estimators (\cite[Rem. 4.2.17]{rieder}), so especially for maximum likelihood estimators, quantiles or least squares estimators. \ \\ \ \\

\section{Compact differentiability of regularized M-functionals}\label{compdiffsec} \

This is the main part of this paper. We will recapitulate the results on asymptotic linearity of unpenalized M-functionals that will be transferred to the regularized case thereafter. \ \\

\subsection{Asymptotic linearity of M-estimators} \

Throughout this section, let $F$ be a distribution on $(\R^p, \B^p)$ and let $X_1,..., X_n \overset{i.id.}{\sim} F$. For some $\Theta \subset \R^p$ ($p$ finite), denote by $\mathcal{C}^p(\Theta)$ the space of all continuous $\R^p-$valued functions on $\Theta$ w.r.t. the supremum norm. A general assumption throughout this paper will be \ \\

\textbf{(A0)} The parametric model $\mathcal{P}=\{P_{\theta} \ | \ \theta \in \Theta\}$ is $L_2-$differentiable and if $\psi_{\theta}$ is an influence curve, it belongs to the set $\Psi_2$. \ \\

Define the function  \begin{center} $ \displaystyle  \eta: \Theta \rightarrow \R^p, \ \ \ \eta(\theta):=\int \varphi(x,\theta) dF(x)=\E_F[\varphi(X,\theta)] $ \end{center} and call its empirical counterpart $Z_n$. The next two assumptions are \ \\

\textbf{(A1)} The parameter space $\Theta \subset \R^p$ is nonempty, compact and equals the topological closure of its interior. \ \\

\textbf{(A2)} The function $\varphi$ satisfies \begin{center} $ \displaystyle  \varphi(x, \cdot) \in \mathcal{C}^p(\Theta) \ F(dx)-\text{a.e.}, \ \ \ \varphi_{\theta}:=\varphi(\cdot, \theta) \in L_2^p(F) \ \forall \theta \in \Theta. $ \end{center}

The following main corollary (\cite[Cor. 1.4.5]{rieder}) makes use of a result of Jain and Marcus (\cite[Thm. 1]{jain75}) which requires  \ \\ 

\textbf{(A3)} There exists a pseudo-distance $d$ on $\Theta$ such that $d(\theta, \theta_0) \rightarrow 0$ as $\theta$ converges to $\theta_0$ and such that the metric integral \vspace{-0.2cm} \begin{center} $ \displaystyle \int_0^1 \sqrt{H(\epsilon)} d\epsilon $ \end{center} is finite. 

\textbf{(A4)} There exists $M \in L_2(F)$ such that \begin{center} $ \displaystyle |\varphi(x,\zeta)-\varphi(x,\theta)| \le d(\zeta, \theta)M(x) \ \forall \zeta, \theta \in \Theta$ \end{center} $F(dx)-$a.e.. \ \\

and of the following theorem \cite[Thm. 1.4.2]{rieder}. 

\begin{thm}[\textbf{Compact differentiability of M-estimators}]\label{riedermain} Under (A1), (A2), assume additionally that: 

\textbf{(A5)} There exists a zero $\theta_0 \in \Theta^{\circ}$ of $\eta$ and $\eta \in \mathcal{C}^p(\Theta)$. Moreover, $\eta$ is locally homeomorphic at $\theta_0$ with bounded and invertible derivative $d \eta(\theta_0)$. \ \\

Then there exists a neighborhood $V \subset \mathcal{C}^p(\Theta)$ of $\eta$ and a functional $T: V \rightarrow \Theta$ satisfying \begin{center} $ \displaystyle f(T(f))=0 \ \forall f \in V. $ \end{center} $T$ is compactly differentiable at $\eta$ with derivative given by \begin{center} $ \displaystyle d_HT(\eta)= -(d\eta(\theta_0))^{-1} \circ \Pi_{\theta_0}, $ \end{center} where $\Pi_{\theta_0}$ is the evaluation functional at $\theta_0$. \end{thm} 

\begin{cor}\label{mestale} Under the assumptions (A0)-(A5), the sequence $(S_n)_n:=(T \circ Z_n)_n$ of M-estimators has the asymptotic linear expansion \begin{center} $  \displaystyle \sqrt{n}(S_n-\theta_0)=\frac{1}{\sqrt{n}} \sum_{i=1}^n \psi_{\theta}(x_i)+o_{(F^n)_*}(n^0)$ \end{center} where the influence function is given by \begin{center} $ \displaystyle  \psi_{\theta}(x):=-(d\eta(\theta_0))^{-1} \varphi(x,\theta_0) $. \end{center} If the $S_n$ are measurable, then asymptotic normality follows, i.e., \begin{center} $ \displaystyle \sqrt{n}(S_n-\theta_0) \circ F^n \overset{w}{\longrightarrow} \mathcal{N}(0,ACA^T) $ \end{center} where \begin{center} $ \displaystyle A:=(d\eta(\theta_0))^{-1}, \ \ \ C:=\E_F[\varphi_{\theta_0} \varphi_{\theta_0}^T]. $ \end{center} \end{cor}

\begin{rem}[\textbf{Fr\'{e}chet differentiability}] \label{hadafre} The proof uses an infinite-dimensional version of the delta method (see e.g. \cite{vawell}, \cite[Thm. 1.3.3]{rieder}) that requires the chain rule. By theorem \ref{compchain}, the functionals have to be at least compactly differentiable. Since the chain rule holds for Fr\'{e}chet differentiable maps, one may ask if the gap between compactly and Fr\'{e}chet differentiable statistical functionals is considerable. The following two examples give an answer. \end{rem}

\begin{ex} \label{quantfre} We refer to \cite[Thm. 1.5.1]{rieder} who shows that for distribution functions $F$ that are continuous in some neighborhood $U$ around $a=F^{-1}(\alpha)$, the location $\alpha-$quantile is compactly but not boundedly differentiable along $\mathcal{C}(U) \cap \mathbb{D}(\R)$, provided that $f(a)>0$, where $\mathbb{D}(\R)$ denotes the Skorohod space, i.e., the space of all real-valued c\`{a}dl\`{a}g functions. \end{ex} 

\begin{ex} \label{wilcoxfre} Another example is given by the functional $T(F,G):=\int FdG$ for distribution functions $F$, $G$. It is shown that this functional is compactly differentiable with Hadamard-derivative \begin{center} $ \displaystyle d_HT(x,y)=\int xdG-\int ydF, $ \end{center} and the empirical version corresponding to the Wilcoxon statistic is compactly differentiable as well (cf. \cite{gill89}, \cite{vawell}). This fact has been used to prove asymptotic linearity of the area under the curve (AUC) and the cross-validated AUC as it has been done in \cite{laan15}. However, \cite{wellner92} showed that $T$ is not Fr\'{e}chet differentiable if one considers the $||\cdot ||_{\infty}-$norm. The given counterexample relies on the fact that in the case of Fr\'{e}chet differentiability, the derivative $d_FT$ coincides with the Hadamard derivative $d_HT$, so $d_HT$ would be the only candidate for $d_FT$, but $d_HT$ does not supply the $o$ term in the first-order expansion in every case. \end{ex}

So, we can summarize that it is reasonable to show the asymptotic linearity by the milder requirement of compact differentiability. \ \\

\subsection{The regression context} \label{vermachlea} \

Fitting a model based on a training set by minimizing some loss function without any restriction generally leads to overfitting, especially in the case of high-dimensional data. This issue has been investigated by Vapnik (\cite{vapstat}) who introduced the structural risk minimization principle which performs the optimization on structures that have finite Vapnik-Chervonenkis dimension. In practice, this idea manifests itself when penalizing the loss function by a regularization term. \ \\

In the regression context, we have a model \begin{center} $ \displaystyle  Y=f(X)+\epsilon $ \end{center} where $(x,y) \in \mathcal{X} \times \mathcal{Y} \subset \R^p \times \R$ for an error term $\epsilon \in \R^n$ with $\epsilon_i \overset{i.id.}{\sim} F_{\epsilon}$ with $\E[\epsilon_i]=0$ and $Var(\epsilon_i)=\sigma^2 \in ]0,\infty[$ for all $i$. The function $f$ might be any measurable function mapping from $\mathcal{X}$ into $\mathcal{Y}$. In this work, we assume that $f$ is an element of the parametric function class \begin{center} $ \displaystyle \mathcal{F}_{\beta}:=\{f_{\beta}(X)=X\beta \ | \ \beta \in \Theta \subset \R^p\}. $ \end{center} 

\begin{rem}[\textbf{Intercept}\index{Model matrix!Intercept}] Note that unless specified otherwise, the first column of the regressor matrix $X$ may only consist of ones, which means that the first component of the parameter is the intercept. \end{rem}

We try to recover the true map $f_{\theta}$ by estimating $\theta$. This is done by defining a loss function $L: (\mathcal{X} \times \mathcal{Y}) \times \Theta \rightarrow [0,\infty[$. For practical applications, we will assume that $L((x,y),f_{\theta})=0$ if $f_{\theta}(x)=y$ as it was done in \cite{christ09}. The penalty term $J_{\lambda}: \Theta \rightarrow [0,\infty[$ that should enforce sparseness of the solution has to satisfy the following conditions: 

\textbf{(A6)} $J_{\lambda}$ is non-negative with $J_{\lambda}(0_p)=0$ and $J_{\lambda}$ is convex. \ \\

The assumption that a regularization term must be non-negative is natural. On the other hand, since it penalizes the model complexity, the assumption that $J_{\lambda}(0_p)=0$ is reasonable since the parameter $\theta=0_p$ leads to a model consisting at most  of the intercept which would not make sense to penalize. The convexity assumption is needed for practical applications to prevent the solution from overfitting and, of course, to guarantee the existence of a unique solution in combination with a convex loss function.  \ \\ 

Then we try to solve  \begin{center} $\displaystyle R(\theta):=\E_F[L((X,Y),f_{\theta})]+J_{\lambda}(\theta)=\int_{\mathcal{X} \times \mathcal{Y}} L((x,y),f_{\theta}) dF(x,y)+J_{\lambda}(\theta)=\min_{\theta \in \Theta}! $ \end{center}  by solving the empirical counterpart of the corresponding Z-equation \begin{center} $ \displaystyle \eta^{\lambda}(\theta):=\int_{\mathcal{X} \times \mathcal{Y}} \varphi((x,y),\theta)dF(x,y)+J_{\lambda}'(\theta)\overset{!}{=}0, $ \end{center} provided that $\partial_{\theta}L=\varphi$ exists and that integration and differentiation can be interchanged. If the penalty term is not of a particular interest, we suppress the superscript and just write $Z_n$ or $\eta$. \ \\

As for $L_2-$differentiability of parametric regression models, we refer to \cite[Thm. 2.4.7]{rieder} for random design and to \cite[Thm. 2.4.2]{rieder} for fixed design of the regressor matrix. \ \\

\subsection{Asymptotic linearity of regularized M-estimators}\label{alesubs} \

For clarity, we write down the following corollary of \ref{mestale} to illustrate that we regard the loss function and the penalty term separately, where the latter will be the one that is more likely to cause problems.  \

\begin{cor}\label{corale} Assume (A0), (A1) and (A3). Let the assumptions (A2) and (A4) be true for \begin{center} $ \displaystyle \varphi((x,y),\theta)+J_{\lambda}'(\theta): \R^p \times \R \times \Theta \rightarrow \R^p $ \end{center}  and let (A5) be true for $\eta^{\lambda}(\theta)$ provided that the derivative exists. Then the asymptotic linear expansion in \ref{mestale} holds with \begin{center} $ \displaystyle \psi(x,y)=-\left(d_H\left(\int \varphi((x,y),\theta)dF(x,y)+J'_{\lambda}(\theta)\right)\bigg|_{\theta=\theta_0}\right)^{-1}\left[\varphi((x,y),\theta)+J_{\lambda}'(\theta)\bigg|_{\theta=\theta_0}\right]. $ \end{center} \end{cor} \ \\

\subsubsection{Compactness assumption of the parameter space} \

We need to justify the assumption of a compact parameter space. This can be done by coercivity arguments. For the following lemma, we refer to \cite{evgrafov} and \cite{levitin}. \ \\ 

\begin{lem} \label{coerc2} Let $f: \mathcal{X} \times \mathcal{Y} \times \Theta \rightarrow \R$ be continuous, where $\mathcal{X} \subset \R^n$, $\mathcal{Y} \subset \R^m$, $\Theta \subset \R^k$. Define $\Xi(x,y):=\argmin_{\theta}(f(x,y,\theta))$. If $f$ is \textbf{coercive w.r.t. $\theta$}, i.e., the sets \begin{center} $ \displaystyle \{\theta \in \Theta \ | \ f(x,y,\theta) \le c\} $ \end{center} are bounded for all $c \in \R$ for every $x \in \mathcal{X}$, $y \in \mathcal{Y}$, then $\min_{\theta}(f(x,y,\theta))>-\infty$ and $\Xi(x,y)$ is nonempty and compact for any $x$, $y$. \end{lem} 

\begin{lem} \label{regcoerc} Let $\mathcal{X}$, $\mathcal{Y}$, $\Theta$ be real vector spaces. Let $L: \mathcal{X} \times \mathcal{Y} \times \Theta \rightarrow [0,\infty[$ be a continuous loss function and let $J_{\lambda}: \Theta \rightarrow [0,\infty[$ be a convex penalty function where $J_{\lambda} \not\equiv 0$. Let $F$ be a distribution on $\B(\mathcal{X} \times \mathcal{Y})$. Then the risk function $\E_F[L((X,Y), \theta)]+J_{\lambda}(\theta)$ is coercive w.r.t. $\theta$, so the parameter space can be restricted to a compact. \end{lem} 

\begin{bew} By convexity, the penalty terms always must satisfy $\lim_{||\theta|| \rightarrow \infty}(J_{\lambda}(\theta))=\infty$ , otherwise it would have to be constantly zero which we excluded by assumption. The coercivity is inherited from the penalty term since the loss function is convex and by linearity of the integral, its expectation is as well, so the risk is coercive w.r.t. $\theta$ (lemma \ref{coerc2}). In fact, we get $R(\theta) \rightarrow \infty$ for $||\theta|| \rightarrow \infty$, so we are allowed to restrict the parameter space to a compact due to lemma \ref{coerc2}. \begin{flushright} $_\Box $ \end{flushright} \end{bew} 

This reasoning is of course not new and has been already done to show the existence of solutions for the Huberized lasso (\cite{lacroix}), for regularized kernel methods in \cite{devito} or in \cite{reyes} for regularized functionals in the context of image restoration. \ \\

It is easy to see that the usual penalty terms like the $l_1-$, $l_2-$ or elastic net penalty are coercive (see p.e. \cite[Cor. 8]{aravkin}). On the other hand, non-convex penalties do not have to be coercive, for example the SCAD penalty (cf. \cite{fanli}) is constant outside a neighborhood of zero whose width depends on the penalty parameter. \ \\ 

In fact, since we are now allowed to assume compactness of the parameter space, we face another potential issue. The compactness assumption leads to the problem that the M-estimator $\hat \theta_n$ may be located at the boundary of $\Theta$. We invoke the idea of one-step estimators from \cite{vaart98} to make the connection with machine learning algorithms. \ \\

Having a $\sqrt{n}-$consistent preliminary solution $\tilde \theta_n$ of the estimating equation $Z_n(\theta)=0$, then an application of the Newton-Raphson algorithm leads to an improved one-step solution \begin{center} $ \displaystyle \hat \theta_n:=\tilde \theta_n-(Z_{n,0}'(\tilde \theta_n))^{-1}Z_n(\tilde \theta_n) $ \end{center} where $Z_0'$ is a regular matrix and $Z_{n,0}'$ is regular and converges in probability to $Z_0'$. The following theorem can be found in \cite[Thm. 5.45]{vaart98}. \ 

\begin{thm} \label{vaartstep} Let the notation be as above. Let the condition that for every constant $M$ it holds that \begin{equation}\label{1step}   \sup_{\sqrt{n}||\theta-\theta_0||<M}\left(||\sqrt{n}(Z_n(\theta)-Z_n(\theta_0))-Z_0'\sqrt{n}(\theta-\theta_0)||\right) \overset{P}{\longrightarrow} 0  \end{equation} be satisfied for a regular matrix $Z_0'$. If it holds additionally that $\sqrt{n}(Z_n(\theta_0))$ converges to some limit, then the One-Step estimator $\hat\theta_n$ is already $\sqrt{n}-$consistent. \end{thm} 

\begin{lem} \label{mqonestep} Let all the notation be as above. Under (A2), (A5) and the additional assumptions 

\textbf{(A7)} The learning procedure is $\sqrt{n}-$consistent, 

\textbf{(A8a)} The function $Z_n$ is twice differentiable w.r.t. $\theta$, 

the One-Step estimator is not located at the boundary of the parameter space. \end{lem}

\begin{bew} Since condition (\ref{1step}) is weaker than differentiability of $Z_n$ at $\theta$, this part is already satisfied by (A8a). The only condition of theorem \ref{vaartstep} that remains to be proven is the convergence of $\sqrt{n}(Z_n(\theta_0))$ to some limit $Z$. But, since we already know by (A7) that the learning algorithm is $\sqrt{n}-$consistent, so $\hat \theta_n-\theta_0=o_F(n^{-1/2})$, hence we get \begin{center} $ \displaystyle \sqrt{n}(\theta_0-\hat \theta_n)=o_F(n^0)$. \end{center} An application of a delta method which is possible under (A8a) provides that \begin{center} $ \displaystyle \sqrt{n}(Z_n(\theta_0)-Z_n(\hat \theta_n))$ \end{center} has a limiting distribution (which is the Dirac measure at zero) and we note that by definition, it holds that $Z_n(\hat \theta_n)=0$, so the convergence of $\sqrt{n}Z_n(\theta_0)$ has been established and theorem \ref{vaartstep} applies. \begin{flushright} $_\Box $ \end{flushright} \end{bew} 

We admit that it is not common to assume learning rates like we did in assumption (A7), but it is more convenient just to assume consistency. Since functions that are too complex may not be able to be approximated with a predetermined rate, this assumption results in the class of approximable functions getting strictly smaller. \ \\

\subsubsection{Twice differentiable Z-function} \

\begin{thm}\label{mqtwicediff} Under the conditions (A0), (A1), (A6), (A7), (A8a) and

\textbf{(A2')} $\varphi^{\lambda}(\cdot, \theta) \in L_2^p(F) \ \forall \theta \in \Theta$, 

\textbf{(A5')} $\eta^{\lambda}(\theta_0)=0$ for a $\theta_0 \in \Theta^{\circ}$ and $\eta^{\lambda}$ is locally homeomorphic at $\theta_0$ with bounded and invertible derivative $d\eta^{\lambda}(\theta_0)$, 

the sequence $(S_n^{\lambda})_n:=(T \circ Z_n^{\lambda})_n$ of regularized M-estimators is asymptotically linear\index{Asymptotic linearity!of regularized M-estimators}\index{Functional differentiability!of regularized M-functionals}. \end{thm} 

\begin{bew} As a byproduct of (A1), we immediately get (A3). This is true since $\R^p$ is a normed space, hence the pseudo-distance is just the standard euclidean norm on $\R^p$ and by boundedness of $\Theta$ and since $p$ is finite, we can conclude that the metric integral is finite. Even more general, the metric integral is finite provided that $d(\theta, \theta_0)=||\theta-\theta_0||_2^{\delta}$ for some $\delta \in ]0,\infty[$ (see \cite[Rem. 1.4.6.b)]{rieder}). \ \\

From twice differentiability of $\eta^{\lambda}$, the first part of (A2) is trivially satisfied and the derivative $\varphi$ of $L$ is continuous w.r.t. $\theta$. By (A6), (A7), (A8a) and lemma \ref{mqonestep}, the assumption that we can restrict the parameter space $\Theta$ to a compact set is justifiable. Using this compactness of $\Theta$, we deduce that the function $\varphi(x,\cdot)$ is Lipschitz-continuous, thus there exists a constant (w.r.t. $\theta$) such that \begin{center} $ \displaystyle |\varphi(x,\xi)-\varphi(x,\theta)| \le L_x||\xi-\theta|| $ \end{center} where the Lipschitz constant $L_x$ must be finite by compactness of $\Theta$. We conclude that (A4) holds. \ \\

 Thus, corollary \ref{mestale} is applicable and we get the desired result.  \begin{flushright} $_\Box $ \end{flushright} \end{bew} \

\subsubsection{Twice continuously differentiable loss function, non-differentiable penalty term} \

If the penalty term is non-differentiable, like the Lasso loss, then we invoke an approximation result of \cite{medina} which uses a maximum theorem of \cite{berge}. This result is exactly what we need in the presence of non-differentiable regularization terms since that despite we cannot assume differentiability, we can at least assume continuity. The following lemma is a combination of \cite[Lemma 2]{medina} and \cite[Prop. 1]{medina}. \  

\begin{lem}[\textbf{Approximating influence curves}]\index{Influence curve/function!Approximation in Sobolev space} \label{medina} Assume that the parameter space $\Theta \subset \R^p$ is compact and that the loss function is twice continuously differentiable w.r.t. $\theta$. If there exists a sequence $(J_{\lambda}^m)_m$ with $J_{\lambda}^m \in C^{\infty}(\Theta)$ that converges to $J_{\lambda}$ in the Sobolev space $W^{2,2}(\Theta)$, i.e., \begin{center} $ \displaystyle ||J_{\lambda}^m-J_{\lambda}||_{W^{2,2}}=\left( \sum_{|\alpha| \le 2} \int_{\Theta} |\partial^{\alpha}(J_{\lambda}^m(\theta)-J_{\lambda}(\theta))|^2 d\theta \right)^{1/2} \longrightarrow 0, $ \end{center} then \begin{center} $ \displaystyle \lim_m(T_m)=T $ \end{center} where $T_m$ denotes the M-functional that intends to find the zero of the Z-equation corresponding to $R_m$ where $R_m$ denotes the risk function where $J_{\lambda}$ is replaced by $J_{\lambda}^m$. This does not depend on the particular choice of the approximating sequence $J_{\lambda}^m$. \end{lem}

\begin{thm} \label{mqapprox} Assume that there exists a sequence $(J_{\lambda}^m)_m$ with $J_{\lambda}^m \in \mathcal{C}^{\infty}(\Omega)$ of regularization functions that converge to $J_{\lambda}$ in the Sobolev space $W^{2,2}(\Theta)$. Under the conditions (A0), (A1), (A2'), (A5'), (A6), (A7) and

\textbf{(A8b)} The loss function $L$ is convex and twice continuously differentiable w.r.t. $\theta$, 

the sequence $(S_n^{\lambda})_n$ of regularized M-estimators is asymptotically linear\index{Asymptotic linearity!of regularized M-estimators}\index{Functional differentiability!of regularized M-functionals}. \end{thm}

\begin{bew} From theorem \ref{mqtwicediff}, we can conclude that the estimator has an asymptotic linear expansion and that it is asymptotically normal if the respective assumptions collected there are satisfied. But since this is just an asymptotic property up to a remainder term of order $n^{-1/2}$, it suffices to approximate $J_{\lambda}$ by $J_{\lambda}^m$ such that the difference in the respective influence functions is negligible, i.e., the difference is already of order $n^{-1/2}$. Note that by continuity of the G\^{a}teaux derivative w.r.t. the direction and by the abovely stated lemma, it holds that $\lim_m(IC(x,T_m,P))=IC(x,T,P)$. \ \\

Finally, we can conclude that we can work with infinitely differentiable penalty terms satisfying the conditions of the previous subsection but that this results in the same asymptotic linear expansion as if we used the true non-differentiable penalty term. Thus, we only need the existence of an approximating sequence of penalty terms. \begin{flushright} $_\Box $ \end{flushright} \end{bew} 

\begin{rem} \cite[Lemma 5.4]{alfons14} showed for the Lasso and a concrete hyperbolic tangent approximation of the penalty term that the influence function of the approximating estimator derived by \cite[Prop. 4.1]{alfons14} converges to the influence function of the Lasso. So, Medina generalized their result with \ref{medina} for any losses and penalties satisfying the given conditions. \end{rem} \

Note again that the main difficulty for non-differentiable regularization terms is the translation of M- to Z-equations. The approximation result elegantly avoids a tedious case-by-case study under which conditions an M-estimator w.r.t. a certain regularized loss function can be written as Z-estimator and provides a universal result. \ \\

\subsection{Extension to ranking} \

In the ranking setting (see for example \cite{clem08}), we assume the same underlying model as in the first part of this section, with the only difference that we have to invoke a joint distribution $F_r: (\mathcal{X} \times \mathcal{Y}) \times (\mathcal{X} \times \mathcal{Y})$ on the measurable space $((\mathcal{X} \times \mathcal{Y}) \times (\mathcal{X} \times \mathcal{Y}), \mathbb{B}((\mathcal{X} \times \mathcal{Y}) \times (\mathcal{X} \times \mathcal{Y})))$ where the notation $F_r$ is introduced to distinguish it from the joint distribution in the previous part of this section. \ \\

In contrast to prediction problems, it is not the goal to recover the true values of the $Y_i$ but just to predict their true order. Therefore, the ranking model can be fitted by defining a ranking loss function $L^{r}: (\mathcal{X} \times \mathcal{Y}) \times (\mathcal{X} \times \mathcal{Y})  \times \Theta \rightarrow [0,\infty[$ which quantifies a ranking loss, that is some loss of a misranking of a pair of instances. Defining a penalty term and the ranking risk analogously to the risk $R$ in section \ref{vermachlea}, the corresponding Z-equation resulting from the problem to minimize the regularized risk is \begin{center} $\displaystyle Z_n^{r,\lambda}(\theta):=\frac{1}{n(n-1)} \mathop{\sum \sum}_{i \ne j} \varphi^r(((X_i,Y_i),(X_j,Y_j)),\theta)+J_{\lambda}'(\theta) \overset{!}{=}0 $\end{center} where $\varphi^r=\partial_{\theta}L^r$ is the score function, hence the first term of $Z_n^{r,\lambda}$ is the empirical counterpart of the expected score.  \ \\ 

We can easily adapt the theorem of \cite{rieder} to the ranking setting and conclude compact differentiability\index{Functional differentiability!of regularized ranking M-functionals|(} of regularized ranking functionals and asymptotic linearity of the sequence $(S_n^{r,\lambda})_n:=(T \circ Z_n^{r,\lambda})_n$ of regularized ranking M-estimators\index{Asymptotic linearity!of regularized ranking M-estimators|(}.  \ \\ 

\begin{thm}\label{mqalerank} Define \begin{center} $ \displaystyle \eta^{r,\lambda}: \Theta \rightarrow \R^p, \ \ \ \eta^{r,\lambda}(\theta):=\int \varphi^r(((x,y),(x',y')),\theta)dF^r(((x,y),(x',y')))+J_{\lambda}'(\theta). $ \end{center}  Then, under the conditions (A0), (A1), (A5), (A6), (A7), (A8a) and \

\textbf{((A2r)')} $\varphi^r(\cdot, \theta)+J_{\lambda}'(\theta) \in L_2^p(F^r) \ \forall \theta \in \Theta$, \

there exists a neighborhood $V \subset \mathcal{C}^p(\Theta)$ of $\eta^{r,\lambda}$ and a functional $T: V \rightarrow \Theta$ satisfying \begin{center} $ \displaystyle f(T(f))=0 \ \forall f \in V. $ \end{center} $T$ is compactly differentiable at $\eta^{r,\lambda}$ with derivative given by \begin{center} $ \displaystyle d_HT(\eta^{r,\lambda})= -(d\eta^{r,\lambda}(\theta_0))^{-1} \circ \Pi_{\theta_0}, $ \end{center} where $\Pi_{\theta_0}$ is the evaluation functional at $\theta_0$. \ \\

Moreover, the sequence $(S_n^{r,\lambda})_n:=(T \circ Z_n^{r,\lambda})_n$ of regularized ranking M-estimators is asymptotically linear\index{Asymptotic linearity!of regularized ranking M-estimators|)}. \end{thm} 

\begin{bew} This directly follows from corollary \ref{mestale} since the optimization is done w.r.t. $\theta$ whereas the dimension of the space of the other arguments of $\varphi$ is not explicitly used in the proof. \begin{flushright} $_\Box $ \end{flushright} \end{bew} 

\begin{rem} It is important to emphasize that ranking loss functions like the hard ranking loss (cf. \cite[Sec. 2]{clem08}) and related losses are not continuous and not convex, so they fail the assumptions of these theorems (however, the hard ranking loss is bounded, so combining it with a suitable regularity term again leads to a coercive target function). Examples for which the regularity conditions hold are smooth convex surrogates of those ranking losses (see \cite{clem13} for an overview). \end{rem}  \ \\

\section{Examples for asymptotically linear estimators in machine learning}\label{exsec} \

The conditions for asymptotic linearity of the regularized M-estimators in the previous section are quite general. The goal of this section is to provide examples for machine learning algorithms to which the derived results can be applied and to specify the required conditions for each procedure. \ \\

\subsection{Lasso} \

Lasso regression (cf. \cite{bu}) is an $l_1-$penalized least squares regression, i.e.,  \begin{center} $ \displaystyle \hat \beta^{lasso}=\argmin_{\beta \in \Theta}\left(\frac{1}{n}||Y-X\beta||_2^2+\lambda ||\beta||_1\right) $. \end{center} The lasso regression results in a shrinkage of the coefficients and in sparsity of the fitted model. The score function for the unregularized loss is given by \begin{equation} \label{scorelasso} \varphi(\cdot,\beta)=\frac{2}{n}X^T(Y-X\beta).  \end{equation} We invoke the approximation of the non-differentiable penalty term. There exists an example of such a smooth penalty term converging to the absolute value in \cite{medina}. \ 

\begin{thm}\label{mqlassoale} Assume (A0), (A1) and

\textbf{((A2$^{Lasso}$)')} The ideal distribution $F$ has finite fourth moments, 

\textbf{(A5')} The true solution $\beta^0$ lies in the interior of $\Theta$ and the derivative $d\eta^{\lambda}(\beta^0)$ is invertible, 

\textbf{(A7$^{Lasso}$)} $||\beta^0||_1=o(\sqrt{n/\ln(p)})$ and that the regularization parameter in dependence of $n$ is chosen in the range of $\lambda_n=\sqrt{\ln(p)/n}$. 

Then the sequence $(S_n^{Lasso})_n:=(T \circ Z_n^{Lasso})_n$ of Lasso estimators is asymptotically linear\index{Lasso!Asymptotic linear expansion}\index{Asymptotic linearity!of the Lasso}. \end{thm} 

\begin{bew} We need to verify the conditions of theorem \ref{mqapprox}. Consider a smooth approximation $J_{\lambda}^m$ of the absolute value in the sense of the Sobolev space $W^{2,2}$, as given in \cite{alfons14} or \cite{medina}, respectively. Then we set \begin{center} $ \displaystyle \tilde J_{\lambda}^m(x):=\sum_i J_{\lambda}^m(x_i) \longrightarrow \sum_i |x_i|=||x||_1 $, \end{center} and thus \begin{center} $ \displaystyle \nabla_x \tilde J_{\lambda}^m(x)=(\partial_{x_1}\tilde J_{\lambda}^m(x),...,\partial_{x_p}\tilde J_{\lambda}^m(x)) \longrightarrow (\sign(x_1),...,\sign(x_p))=\nabla_x||x||_1 $ \end{center} and the (component-wise) convergence of the Hessian holds as well due to the properties of $W^{2,2}$. For this idea, we refer to \cite[Lemma 5.4]{alfons14}. The loss function and the approximating penalty term are smooth, hence (A8b) is satisfied and lemma \ref{medina} is applicable. \ \\

The target function is coercive w.r.t. $\beta$ (see \ref{coerc2}). This holds because as $||\beta|| \rightarrow \infty$, the penalty will tend to infinity and so does the target function. Note that this does not hold for the loss function itself since $||\beta|| \rightarrow \infty$ can result in a small loss. One may argue that even in the penalized case, it can happen that $||(x,y,\beta)|| \rightarrow \infty$ without resulting in the target function growing as well. If for example $Y=0$ and $X$ is large, then $Y=X\beta$ for $\beta=0_p$. But in this case, we do not lose anything if we restrict the parameter space. Furthermore, we can write the optimization problem in the form \begin{center} $ \displaystyle \min(||Y-X\beta||_2^2/n) \ \ \ s.t. \ \ \ ||\beta||_1 \le c_{\lambda} $ \end{center} for some constant $c_{\lambda}$ depending on $\lambda$. So we have a convex optimization problem with a continuous, strictly convex and coercive target function, so by \cite{werner06}, there definitely exists a solution $\beta^0$ of $\eta^{\lambda}$ and the local homeomorphicity around the solution follows. \ \\

Combining ((A2$^{Lasso}$)') with equation (\ref{scorelasso}), we derive that the score function is square-integrable w.r.t. the distribution $F$. Then (A2') is satisfied and by boundedness of the integral by the previous assumption and by compactness of the parameter space, this derivative is bounded. \ \\ 

The Lasso is generally inconsistent, but under (A7$^{Lasso}$), it follows from \cite{bu} that the Lasso is $\sqrt{n}-$consistent in this case. Note that despite we solve a convex optimization problem assuming that the true solution is already located in the interior of $\Theta$, that does not suffice to guarantee that the computed solution does not lie on the boundary of the parameter space. Finally, theorem \ref{mqapprox} applies and the assertion is proven. \begin{flushright} $_\Box $ \end{flushright} \end{bew} \ 

\subsubsection{Elastic net} \

The elastic net (cf. \cite{zou05}) can be regarded as a compromise between Lasso and Ridge regression. Given two penalty parameters $\lambda_1, \lambda_2$, the elastic net solution is given by  \begin{center} $ \displaystyle \hat \beta^{EN}=\argmin_{\beta}\left(\frac{1}{n} ||Y-X\beta||_2^2+\lambda_1 ||\beta_1||_1+\lambda_2 ||\beta||_2^2 \right) $ \end{center} and by defining $\alpha:=\frac{\lambda_2}{\lambda_1+\lambda_2}$, this can be rewritten as a convex combination of $l_1-$ and $l_2-$ penalties, i.e., \begin{center} $ \displaystyle \hat \beta^{EN}=\argmin_{\beta}\left(\frac{1}{n} ||Y-X\beta||_2^2+(1-\alpha) ||\beta_1||_1+\alpha ||\beta||_2^2 \right) $ \end{center} where $(1-\alpha)||\beta||_1+\alpha ||\beta||_2^2$ is referred to as the elastic net penalty. \ 

\begin{cor} \label{mqenorthale} Under the assumptions of theorem \ref{mqlassoale}, the sequence $(S_n^{orthEN})_n:=(T \circ Z_n^{orthEN})_n$ of elastic net estimators with orthonormal design is asymptotically linear\index{Asymptotic linearity!of the Elastic net}. \end{cor}

\begin{bew} Note that for orthonormal design, the EN solution is just a rescaled Lasso solution with factor $\frac{1}{1+\lambda_2}$. In this case, we can simply rescale the influence function derived in \cite{alfons14}, proving the result. \begin{flushright} $_\Box $ \end{flushright} \end{bew} 

\begin{cor} \label{mqenale} Under the assumptions of theorem \ref{mqlassoale}, the sequence $(S_n^{EN})_n:=(T \circ Z_n^{EN})_n$ of elastic net estimators is asymptotically linear\index{Asymptotic linearity!of the Elastic net}. \end{cor}

\begin{bew} It is shown in \cite{zou05} that the elastic net can be rewritten as a special Lasso with the augmented data \begin{center} $ \displaystyle X^*:=\frac{1}{\sqrt{1+\lambda_2}}\begin{pmatrix} X \\ \sqrt{\lambda_2}I_p \end{pmatrix}, \ \ \ y^*:=\begin{pmatrix} y \\ 0_p \end{pmatrix} $ \end{center} and the penalty $\gamma:=\frac{\lambda_1}{\sqrt{1+\lambda_2}}$. If $\hat \beta^{Lasso}$ is the respective Lasso solution, the elastic net solution is a rescaling with factor $\frac{1}{1+\lambda_2}$ as before. \ \\ 

Using these results and the idea of \cite{alfons14}, the respective influence curve of the elastic net can been computed by just adapting the already calculated influence curve and theorem \ref{mqlassoale}. \begin{flushright} $_\Box $ \end{flushright} \end{bew}

\subsubsection{Adaptive Lasso} \

The adaptive Lasso (cf. \cite{zou06}) is a two-stage estimator that first computes the standard Lasso estimator, denoted by $\hat \beta^{init}$, and then in a second step, one minimizes \begin{center} $ \displaystyle \frac{1}{n} ||Y-X\beta||_2^2+\lambda \sum_j \frac{|\beta_j|}{|\hat \beta^{init}_j|}. $ \end{center} 

Borrowing the consistency requirements for the adaptive Lasso from \cite{zou06}, we have the following result. \ 

\begin{thm}\label{mqalassoale}  Assume (A0), (A1) and

\textbf{((A2$^{Lasso}$)')} The ideal distribution $F$ has finite fourth moments, 

\textbf{(A5')} The true solution $\beta^0$ lies in the interior of $\Theta$ and the derivative $d\eta^{\lambda}(\beta^0)$ is invertible,

\textbf{(A7$^{ALasso}$)} The regularization parameter in dependence of $n$ satisfies $\lambda_n=o(\sqrt{n})$ and $\lambda_n n^{(\gamma-1)/2} \rightarrow \infty$ for $\gamma>0$. 

Then the sequence $(S_n^{adapt})_n:=(T \circ Z_n^{adapt})_n$ of Adaptive Lasso estimators is asymptotically linear\index{Asymptotic linearity!of the Adaptive Lasso} and the influence curve\index{Influence curve/function!of the Adaptive Lasso} of the $j-$th component of $\hat \beta^{adapt}$ is given by \begin{center} $ \displaystyle IC((x_0,y_0), \hat \beta_j^{adapt}(\lambda), P_{\theta_0})=\begin{cases} 0, \ \ \ \hat \beta_j^{init}(\lambda)=0 \\ 0, \ \ \ \hat \beta_j^{adapt}(\lambda)=0 \\ IC((x_0,y_0), \hat \beta_j^{Lasso}(\lambda/|\hat \beta_j^{init}(\lambda)|)), \ \ \ \text{otherwise} \end{cases}. $ \end{center} where we denote by $\hat \beta^{Lasso}(\lambda)$ the Lasso estimator using the penalty factor $\lambda$. \end{thm} 

\begin{bew} Obviously, if $\hat \beta^{init}_j=0$, we immediately know that $\hat \beta^{adapt}_j=0$. Hence, if we have the initial solution, we can rewrite the adaptive Lasso optimization problem as \begin{center} $ \displaystyle \hat \beta^{adapt}_{\hat S_{\lambda}^{init}}=\argmin_{\beta_{\hat S^{init}(\lambda)}} \left(\frac{1}{n} \sum_{i=1}^n \sum_{j \in \hat S^{init}(\lambda)} (Y_i-X_{ij}^T \beta_j)^2+\lambda \sum_{j \in \hat S^{init}(\lambda)} \frac{|\beta_j|}{|\hat \beta_j^{init}|} \right)$ \end{center} where $\hat S^{init}(\lambda):=\{j \ | \ \hat \beta^{init}_j \ne 0\}$. Then this optimization problem is just a Lasso optimization problem with a weighted penalty term which can be approximated coordinate-wisely in the spirit of \cite{medina}. \ \\

The corresponding influence function has implicitly been derived in \cite{medina}. Note that by our method, we would only derive $|\hat S^{init}(\lambda)|$ components of the influence function. However, it was proven in \cite{alfons14} that the components of the influence function corresponding to the coefficients that are excluded from the model are zero, i.e., if the Lasso in the first step already sets some coefficients to zero, the final coefficients will be zero, so we can just plug in zeroes into the respective components of the influence function, providing the usual asymptotic linear expansion. \ \\

Since the first step does not compute the final non-zero coefficients but just regularizing weights, its influence implicitly arises in this expansion as a factor, leading to the stated result. \begin{flushright} $_\Box $ \end{flushright} \end{bew} 

\begin{rem}[\textbf{Partial influence curves}]\index{Influence curve/function!Partial} \label{alassorem} Note that the influence curves derived in theorem \ref{mqalassoale} correspond to the concept of ''partial'' influence functions (see definition \ref{alepic}). This is true since in the proof of theorem \ref{mqalassoale}, we are implicitly using the smooth transformation $\beta \mapsto \beta_{\hat S_{\lambda}^{init}}$ to derive the components of the influence curve corresponding to the coefficients that not already have been excluded from the model in the initialization step. In other words (after suitable renumeration of the columns), we get the matrix $D_{\beta_{\hat S^{init}(\lambda)}}:=(\text{diag}(1,\hat s^{init}),0_{p-\hat s^{init}}) \in \R^{\hat s^{init} \times p}$ where $\hat s^{init}:=|\hat S^{init}(\lambda)|$. \end{rem}

\begin{rem}[\textbf{Asymptotic normality}] Note that the additional assumption of measurability of the sequences of estimators provides asymptotic normality of the estimating sequence due to corollary \ref{mestale}. Of course, we are not the first ones with results on asymptotic normality. See for example \cite[Sec. 3]{loh17} where Corollary 1 shows under which conditions regularized M-estimators of a very general form, including the Lasso, are asymptotically normal. \end{rem} \ \\

\section{Data-driven penalty parameters} \

It is common that asymptotic results for regularized methods allow for the case that the regularization parameter is data-driven which manifests itself in a sequence $(\lambda_n)_n$ of regulariztion parameters. The same is true for Boosting where the amount of regularization does not depend on a penalty parameter but implicitly on the number of iterations such that a diverging sequence of iterations is the analogue to a sequence $(\lambda_n)_n$ with $\lambda_n \rightarrow 0$ for $n \rightarrow \infty$. \ \\

To keep the asymptotic results valid uniformly for $n \rightarrow \infty$, results for Lasso methods as in \cite{bu} or \cite{zou06} and for Boosting methods as for example in \cite{bu06} require penalty parameters which fall into a suitable range in dependence of $n$ or numbers of iterations that grow sufficiently slowly w.r.t. $n$. \ \\

As for our results, we would need a suitable degree of approximation, i.e., a suitable sequence $(m_n)_n$, leading to a sequence of regularization terms of the form $(J_{\lambda_n}^{m_n})_n$, to get similar statements. \ \\

In fact, we already used sequences $(m_n)_n$ implicitly when proving theorem \ref{mqapprox}. Our argument was to set $m_n$ sufficiently large to get a degree of approximation that leads to an error term which is already absorbed by $o_P(n^{-1/2})$. \ \\

If we are concerned about sequences of penalty terms, we essentially need to have a sequence $(m_n)_n$ which again grows sufficiently fast to keep the error term small enough. Since we assumed that $J_{\lambda}$ is approximable by a sequence of smooth penalty terms $J_{\lambda}^m$ and since $\lambda$ usually just enters as a factor, a diminishing sequence of regularization terms still keeps the approximability valid since for smaller penalty parameters, the regularization term gets ,,less wiggly'', so we assume that for fixed $n$, one would generally need a smaller number $m$ for a smaller $\lambda$ than for a large $\lambda$. \ \\

A general approximation to the best of our knowledge is out of reach, however, for a given penalty term with a given approximation sequence, one could derive conditions for the sequence $(m_n)_n$ according to the sequence of regularization terms. In particular, this holds for the Lasso and the Adaptive Lasso.                        \ \\  \ \\

\section{Conclusion} \

We studied the conditions under which a regularized M-estimator can be asymptotically linearized. We provide a general theory for the asymptotically linear expansion of such estimators and gave concrete examples of machine learning algorithms which our theory includes. Of course, from the asymptotic linear expansion, the asymptotic normality can be directly derived. \ \\

Influence curves for a wide range of estimators have already been derived in the literature, but it does not suffice just to write down the asymptotic linear expansion without having proved the necessary regularity conditions of the corresponding statistical functional. This means that the potential linearization of the estimators in related work was not theoretically founded without the rigorous theory that we provided. \ \\

However, we concentrated on linear regression models in this work. An extension to other areas of machine learning will be a subject of future work. \ \\

\section*{Acknowledgements} \

The results presented in this paper are part of the author's PhD thesis (\cite{TWphd}) supervised by P. Ruckdeschel at Carl von Ossietzky University Oldenburg. \ \\  \ \\

\renewcommand\refname{References}
\bibliography{Biblio}
\bibliographystyle{abbrvnat}
\setcitestyle{authoryear,open={((},close={))}}

\end{small}

\end{document}